\documentclass[12pt]{article}
\usepackage{amsfonts, amsmath, amssymb}

\newtheorem{thm}{Theorem}

\newtheorem{defi}{Definition}

\def\C{\mathbb C}

\def\P{\mathbb P}

\def\R{\mathbb R}
\def\Z{\mathbb Z}
\def\O{\mathcal O}

\def\defq{\stackrel{def}{=}}

\def\ra{\rightarrow }
\def\hX{\widehat{X}}

\def\noi{\noindent}

\title{Blowing-up points on l.c.K. manifolds.}
\author{Victor Vuletescu}
\date{}

\begin{document}
\maketitle

\begin{abstract}
It is a classical result, due to F. Tricceri, that the blow-up of a manifold of
locally conformally K\"ahler (l.c.K. for short) type at some point is again of
l.c.K. type. However, the proof given in \cite{Tric} is somehow unclear. We give a
different argument to prove the result, using "standard tricks" in algebraic
geometry.
\end{abstract}

\noi {\bf Keywords:}
Blow-up of a manifold at a point, locally conformally K\"ahler manifold, Lee  form.

\noi {\bf 2000 Mathematics subject classification:} 53C55, 14E99
\section{Introduction}

We begin by recalling the basic definitions and facts; details can be found for
instance in the book \cite{DO}.

\begin{defi}Let $(X, J)$ be a complex manifold. A hermitian metric $g$ on it 
is called {\em locally conformally K\"ahler}, l.c. K. for short, if there exists
some open cover ${\mathcal U}=\{U_\alpha\}_{\alpha \in A}$ of $X$  such that for
each $\alpha \in A$ there is some smooth function $f_\alpha$ defined on $U_\alpha$
such that the metric $g_\alpha=e^{-f_\alpha}g$ is K\"ahler.

A complex manifold $(X, J)$ will be called of l.c.K. type if it admits an l.c.K. metric
\end{defi}

Letting $\omega$ to be the {\em K\"ahler form} associated to $g$ by $\omega(X,
Y)=g(X, JY),$
one can immediately show that the above definition is equivalent to the existence of
a closed $1-$form $\theta$ such that $d\omega=\theta\wedge \omega.$ The form
$\theta$ is called {\em the Lee form} of the metric $g.$
It is almost immediate to see that $\theta$ is closed; iy is exact iff the metric $g$ is global conformally equivalent to a  K\"ahler metric. Usually, by an l.c.K. manifold one understands a hermitian manifold whose metric is not globally conformally K\"ahler. 
In particular, the first Betti number of an l.c.K. manifold is always strictly positive; more, for  compact  Vaisman manifolds (l.c.K. with parallel Lee form) the fundamental group fits into an exact sequence 
\[
0 \rightarrow G \rightarrow \pi_1(M) \rightarrow \pi_1(X) \rightarrow 0
\]
where $\pi_1(X)$ is a fundamental group of a K\"ahler  orbifold,
and $G$ a quotient of $\Z^2$ by a subgroup of rank $\geq 1$ (see \cite{ov}). 
Moreover, the l.c.K. class is not stable to small deformations: some Inoue surfaces do not admit l.c.K. structures and they are complex deformations of other Inoue surfaces with l.c.K. metrics (see \cite{Tric}, \cite{bel}).

However, l.c.K. manifolds share with the K\"ahler ones the property of being closed under   blowing-up points. To can state the result, let  
$X$ be a a complex manifold and  $P\in X$ some  fixed point. We denote by 
$\hX$ the manifold obtained by blowing-up $P,$ by $c:\hX\ra X$ the blowing-up map and
$E$ the exceptional divisor of $\pi$ (i.e. $E=c^{-1}(\{P\})$). The goal is to prove
the following

\begin{thm}If the complex manifold $X$ carries an l.c.K. metric, then so does its blow-up $\hX$ at any point.

\end{thm}

The result was stated in  \cite{Tric}, but the proof in this paper has a gap. 

For the sake of completeness, we include in the next section some basic facts about
blow-up's of points on complex manifolds. Eventually, in the last section we prove
the theorem.
%\newpage
\section{Basic facts about blow-up's of points.}

This section is entirely standard and is almost an verbatim reproduction of facts
from classical texts, as for instance \cite{GrHa}.

Let $X$ be a complex, $n-$dimensional manifold. Let $P\in X$ be a point; choose a
holomorphic local coordinate system $(x_1,\dots, x_n)$ defined in some open
neighborhood $U$ of  $P$ such that $x_1(P)=\dots=x_n(P)=0.$ Consider the manifold
$U\times \P^{n-1}(\C)$ and assume $[y_1:\dots:y_n]$ is some fixed homogenous
coordinate system on $\P^{n-1}(\C).$ Let $\widehat{U}\subset U\times \P^{n-1}(\C)$
be the closed subset defined by the system of equations $x_iy_j=x_jy_i, 1\leq
i<j\leq n.$ One can check that $\widehat{U}$ is actually a submanifold of $U\times
\P^{n-1}(\C).$ Moreover, the restriction of the projection onto the first factor
$c:\widehat{U}\ra U$ has the following properties: the fiber of $c$ above $P$,
$c^{-1}(\{P\}),$ is a submanifold $E$ of $\widehat{U}$   which is biholomorphic to
$\P^{n-1}(\C)$ and the restriction of $c$ at $\widehat{U}\setminus E$ defines a
biholomorphism between $\widehat{U}\setminus E$ and $U\setminus \{P\}.$ Using it, we
can glue $\widehat{U}$ to $X$ along $U\setminus \{P\}.$ The resulting manifold will
usually be denoted by $\hX;$ the map $c$ above extends obviously to a map - denoted
by the same letter- $c:\hX\ra X.$ Notice that on one hand $c$ is a biholomorphic map
between $\hX\setminus E$ and $X\setminus \{P\}$ and, on the other hand, $c$
"contracts" $E$, i.e. $c(E)=\{P\}$ ($E$ is called accordingly the "exceptional
divisor" of $c$).

Let now $y\in \hX$ be some point. If $y\not\in E$, then the tangent map $$c_{*,
y}:T_y(\hX)\ra T_{c(y)}(X)$$ is a isomorphism, while if $y\in E$ then the rank of
this map is one and  its kernel consists of those vectors that are tangent at $y$ to
$E, $ i.e. $Ker(c_{*, y})=T_y(E).$

Next, recall that to each closed complex submanifold $E$ of codimension one of some
complex manifold $X$ one can associate a holomorphic vector bundle, usually denoted
$\O_X(E);$ see e.g. \cite{GrHa}, Chapter 1, Section 1. If one chooses a hermitian
metric $h$ in $\O_X(E)$ there exists and is unique a linear connection $D$ in the
vector bundle which is also compatible with the complex structure (see e.g. the
Lemma  on page 73, \cite{GrHa}). The  curvature $\Omega_E$ of this connection is a
closed $(1,1)-$form.

We shall next exemplify the computation of the curvature of a metric connection in 
the special case we are interested in, namely when $E$ is the exceptional divisor of
some blow-up. So let $X$ be a manifold, $P\in X$, $U$ a coordinate neighborhood of
$P$ as in the beginning of the section and  $\hX$ the blow-up of $X$ at $P.$ For
$\varepsilon$ small enough set 
$$U_{2\varepsilon}\defq Q\in U\vert \; \vert x_i(Q)\vert<2e\text{ for all }
i=1,\dots, n\}.$$

Let $\pi':U\times \P^{n-1}(\C)\ra \P^{n-1}(\C)$ be the projection onto the the
second factor; then $\O_{\widehat{U}}(E)=\pi'^*(\O_{\P^{n-1}(C)}(-1)).$
Let $\omega_{FS}$ be the K\"ahler form of the Fubini-Study metric on $\P^{n-1}(\C);$
then $-\omega_{FS}$ is the curvature of the canonical connection of the natural
metric $h$ in the tautological line bundle $\O_{\P^{n-1}(C)}(-1).$
Let $h'\defq \pi'^*(h)$ be the induced metric in $\O_{\widehat{U}}(E);$ then its
curvature will be $\pi'^*(-\omega_{FS}).$
On the other hand, the line bundle $\O_{\hX}(E)$ is trivial outside $E;$ fix a
nowhere vanishing section $\sigma$ of it and let $h"$ be the unique metric making
$\sigma$ into a unitary basis. Let now $\varrho_1, \varrho_2$ be a partition of
unity such that $\varrho_1\equiv 1$ on $U_\varepsilon$ and $\varrho_1\equiv 0$
outside $U_{2\varepsilon}$   and respectively $\varrho_2\equiv 0$ on $U_\varepsilon$
and $\equiv 1$ outside $U_{2\varepsilon}.$
Let $h=\varrho_1 h'+\varrho_2h";$ it is a hermitian metric on $\O_X(E).$ Its
curvature will be zero outside $U_{2\varepsilon}$ since $h=h"$ there. In
$U_\varepsilon$, its curvature will be the pull-back (via $\pi'$) of $-\omega_{FS}$,
hence it is semi-negative definite; moreover, its restriction to $E$ will be
negative definite on vectors that are tangent along $E$, since the restriction of
$\pi'$ to $E$ is a biholomorphism between $E$ and $\P^{n-1}(\C).$.

%\newpage
\section{Proof of the theorem.}

{\bf Proof.}
First, let us fix the terminology. We will say that a $(1,1)-$form $\omega$ on a
complex manifold $(M, J_M)$ is positive (semi-)definite if for any point $m\in M$
and any non-zero tangent vector $v\in T_mM$ one has $\omega(v, J_Mv)>0$
(respectively $\geq 0$), in other words if it is the K\"ahler form of some hermitian
metric on $M.$

Let now $\omega$ be the K\"ahler form of an l.c.K. metric on $X.$ We see
$c^*(\omega)$  is a  $(1,1)-$ form on $\hX$ which is positive definite on
$X\setminus E$ and satisfies $dc^*(\omega)=c^*(\theta)\wedge c^*(\omega)$, where
$\theta$ is the Lee form of the given l.c.K. metric on $X.$ As $E$ is simply
connected we see  (e.g. by using Lemma 4.4 in \cite{Tric}) there exists an open
neighborhood $U$ of $E$ and a smooth function $f:\hX\ra \R$ such that $\omega'\defq
e^fc^*(\omega)$ satisfies $d\omega'=\theta'\wedge \omega'$ and such that
$\theta'_{\vert U}\equiv 0.$

On the other hand, we can find a hermitian metric  in the holomorphic line bundle
$\O_{\hX}(E)$ on $\hX$ associated to $E$ such that the  curvature  $\Omega_E$ of its
canonical connection  is negative  definite along $E$
(i.e. $\Omega_E(v, J_{\hX}v)<0$ for every non-vanishing vector  $v\in T_P(E)$ and
for every $P\in E$), is negatively semidefinite at points of $E$ (i.e. $\Omega_E(v,
J_{\hX}v)\leq 0$ for any $P\in E$ and any $v\in T_P(\hX)$) and is zero outside $U$
(cf. e.g. \cite{GrHa}, pp 185-187). Notice that $\Omega_E$ is closed.

We infer that for some positive integer $N$ the $(1,1)-$form
$h\defq N\omega'-\Omega_E$ is positive definite.

Indeed, this is obvious outside $U$ as $\Omega_E$ vanishes here and $N\omega'$ is
positive definite for any $N>0$.

Along $E,$ as both $\omega'$ and $-\Omega_E$ are positive semidefinite, we have only
to check the definiteness of $h.$ Let $y\in E$ be some point and $v\in T_y(\hX).$ 
Assume $h(v, J_{\hX}v)=0;$ since both $\omega'$ and $-\Omega_E$ are positive
semidefinite, we get $\omega'(v, J_{\hX}v)=\Omega_E(v, J_{\hX}v)=0.$ But $\omega'(v,
J_{\hX}v)=0$ implies  $c^*(\omega)(v, J_{\hX}v)=0;$  so $\omega(c_{*, y}(v),
J_{\hX}c_{*, y}(v))=0$ hence $v\in Ker(c_{*, y}).$ As $Ker(c_{*, y})=T_y(E)$ we get
that $v\in T_y(E);$ but as $-\Omega_E(v, J_{\hX}v)=0$ we see that $v=0$

To check the assertion on $U$, it suffices to notice that for each point $x$ in $U$
there exists some $n_x$ such that $N\omega'-\Omega_E$ is positive definite at $x$
for all $N\geq n_x$, hence also positive definite on a neighborhood of $x;$ since
$U$ is relatively compact, we can cover it by finitely many such neighborhoods, and
take the maximum of the corresponding $n_x'$s.

Last, let us see that $N\omega'-\Omega_E$ is l.c.k. One has
$$d(N\omega'-\Omega_E)=Nd\omega'=\theta'\wedge N\omega'.$$
But $\theta'\wedge \Omega_E=0$ since their  supports are disjoint, so we  have
$$d(N\omega'-\Omega_E)=\theta'\wedge N\omega'-\theta'\wedge \Omega_E
=\theta'\wedge(N\omega'-\Omega_E).$$

\medskip

\noi {\bf Acknowledgments.} I wish to thank L.Ornea  and I. Vaisman for useful
discussions; also, I'm especially grateful to V. Br\^{\i}nz\u anescu for a careful
reading of a preliminary version of this paper.

\begin{flushright}
Universitatea Bucure\c sti,

Facultatea de Matematic\u a \c si Informatic\u a,

Str. Academiei 14,

010014, Bucure\c sti, Rom\^ania.

Email: vuli@fmi.unibuc.ro

\end{flushright}

\end{document}